\documentclass[12pt]{article}
\usepackage{amssymb, amsmath, latexsym, a4}
\usepackage[latin1]{inputenc}
\usepackage[all]{xy}

\begin{document}
\newcommand{\A}{\mathrm A}

\newcommand{\qed}{\hfill $\Box $}
\def \G{\mathop{\Gamma }\nolimits}
\def \Y{\mathop{\cal Y }\nolimits}
\def \n{\mathop{\underline n }\nolimits}
\def \m{\mathop{\underline m }\nolimits}
\newcommand{\Ga}{\raisebox{0.1mm}{$\Gamma$}}
\def\S{\smallskip \par}
\def\M{\medskip \par}
\def\ll{ \lambda}
\def\B{\bigskip \par}
\def\BB{\bigskip \bigskip \par}
\newcommand{\Q}{\mathbb{Q}}
\newcommand{\F}{\mathbb{F}}

\newcommand{\Z}{\mathbb{Z}}
\def\t{\otimes }
\def\ee{\epsilon }
\newcommand{\pn}{\par \noindent}
\def \H{\mathop{\sf H}\nolimits}
\def \Im{\mathop{\sf Im}\nolimits}
\def \Ker{\mathop{\sf Ker}\nolimits}
\def \Coker{\mathop{\sf Coker}\nolimits}
\def \Hom{\mathop{\sf Hom}\nolimits}
\def \Ext{\mathop{\sf Ext}\nolimits}
\def \Tor{\mathop{\sf Tor}\nolimits}
\newcommand{\tri}{ H_{triple}}
\newtheorem{De}{Definition}[section]
\newtheorem{Th}[De]{Theorem}
\newtheorem{Pro}[De]{Proposition}
\newtheorem{rem}[De]{Remark}
\newtheorem{Le}[De]{Lemma}
\newtheorem{Co}[De]{Corollary}
\newtheorem{Rem}[De]{Remark}
\newtheorem{Ex}[De]{Example}
\newcommand{\Def}[1]{\begin{De}#1\end{De}}
\newcommand{\Thm}[1]{\begin{Th}#1\end{Th}}
\newcommand{\Prop}[1]{\begin{Pro}#1\end{Pro}}

\newcommand{\Lem}[1]{\begin{Le}#1\end{Le}}
\newcommand{\Cor}[1]{\begin{Co}#1\end{Co}}
\newcommand{\Rek}[1]{\begin{Rem}#1\end{Rem}}

\newcommand{\ele}{\cal L} \newcommand{\re}{\cal R}  \newcommand{\pe}{\cal P}  \newcommand{\ene}{\cal N}
 \newcommand{\de}{\cal D}  \newcommand{\os}{\cal O}   \newcommand{\lr}{\rm Lie--Rinehart}

\bigskip
\bigskip

\bigskip
\centerline {\large {\bf {Triple Cohomology of Lie--Rinehart Algebras}}}
\centerline {\large {\bf {and the Canonical  Class of  Associative Algebras}}}
\bigskip

\centerline{ \bf J.M. Casas$^1$, M. Ladra$^2$ and T. Pirashvili$^3$}

\bigskip
\centerline{(1) Dpto.  Matem\'atica Aplicada I, Univ. de Vigo, 36005 Pontevedra, Spain}

\smallskip
\centerline{(2) Dpto.  \'Algebra, Univ. de Santiago, 15782 Santiago de Compostela, Spain}
\smallskip
\centerline{(3) Math. Inst., Alexidze str.1, Tbilisi, 0193, Republic of Georgia}

\bigskip
\date{}

\bigskip

\bigskip

\par
{\bf Abstract} {\footnotesize We introduce a bicomplex which computes the
triple cohomology of Lie--Rinehart algebras. We prove that the triple
cohomology is isomorphic to the Rinehart cohomology \cite{Ri} provided the
Lie--Rinehart algebra is projective over the corresponding commutative algebra.
As an application we construct a canonical class in the third dimensional
cohomology corresponding to an associative algebra. }


{\it A. M. S. Subject Class. (2000): 18G60, 16W25, 17A99.}

\section{Introduction}
Let A be a commutative algebra over a field $K$. A Lie--Rinehart algebra is a
Lie $K$-algebra, which is also an A-module and these two structures are related
in an appropriate way \cite{Hu}. The leading example of Lie--Rinehart algebras is
the set Der(A) of all $K$-derivations of A. Lie-Rinehart algebras are algebraic
counterpart of Lie algebroids \cite{mak}.

The cohomology  $H^*_{Rin}(\ele,M)$ of a Lie--Rinehart algebra $\ele$ with
coefficient in a Lie-Rinehart module $M$ first was defined by Rinehart
\cite{Ri} and then was further developed by Huebschmann \cite{Hu}. However
these groups have good properties only in the case, when $\ele$ is projective
over A.  In this paper following to \cite{pira} we introduce a bicomplexes
$C^{**}( {\mathrm A}, \ele, M)$, whose cohomology
 $H^*({\mathrm A},\ele,M)$ is isomorphic to $H^*_{Rin}(\ele,M)$ provided
 $\ele$ is projective as an A-module. It turns out, that for general $\ele$ the
 group $H^*( {\mathrm A},\ele,M)$ is isomorphic to the triple cohomology
 of Barr-Beck \cite{bb} applied to Lie--Rinehart algebras. We also prove
 that for general $\ele$, unlike to the Rinehart cohomology
 $H^*_{Rin}({\ele},M)$, the groups $H^*( {\mathrm A},\ele,M)$
 in dimensions two and three  classify all abelian
 and crossed extensions of $\ele$ by $M$.

As an application we consider the following situation. Let $S$ be an
associative algebra over a field $K$. We let $H^*(S,S)$ be the Hochschild
cohomology of $S$. It is well-known that $H^1(S,S)$ is a Lie $K$-algebra. It
turns out that $H^1(S,S)$ is in fact a Lie-Rinehart algebra over A, where
A=$H^0(S,S)$ is the center of $S$. Thus one can consider the cohomology $H^*(
{\mathrm A}, H^1(S,S),{\mathrm A})$. We construct an element $$o(S)\in
H^3({\mathrm A}, H^1(S,S),{\mathrm A})$$ which we call  {\it the canonical
class of $S$}. We prove that $o(S)$ is a Morita invariant. The construction of
$o(S)$ uses crossed modules of Lie-Rinehart algebras introduced in
\cite{clpcr}.

\section{Preliminaries on Lie-Rinehart algebras}
The material of this section is well-known. We included it in order to fix
terminology, notations  and main examples. In what follows we fix a field $K$.
All vector spaces are considered over $K$. We write $\t$ and $\Hom$ instead of
$\t_K$ and $\Hom_K$.

\subsection{Definitions, Examples}
 Let A be a commutative algebra over
a field $K$. Then the set  Der(A) of all $K$-derivations of A is a Lie
$K$-algebra and an A-module simultaneously. These two structures are related by
the  following identity
$$[D,aD'] = a[D,D'] +D(a)D', \quad  D, D' \in {\rm Der(A)}.$$ This leads to the
following notion, which goes back to Herz under the name "pseudo-alg\`ebre de
Lie" (see \cite{He}) and which are algebraic counterpart of Lie algebroids
\cite{mak}.
\begin{De}{\rm A {\it Lie-Rinehart algebra} over A consists with a Lie $K$-algebra
${\ele}$ together with an A-module structure on ${\ele}$ and a map
$$\alpha : {\ele} \to {\rm Der(A)}$$ which is simultaneously a Lie algebra
 and A-module  homomorphism such that $$[X,aY] = a[X,Y]+X(a)Y.$$ Here
$X,Y \in {\ele}$, $a \in {\rm A}$ and we  write $X(a)$ for $\alpha(X)(a)$
\cite{Hu}. These objects are also known as $(K, {\mathrm A})$-Lie algebras
\cite{Ri} and $d$-Lie rings \cite{pal}. }
\end{De}
 Thus  Der(A) with $\alpha={\sf Id}_{\mathrm Der(A)}$ is a  
Lie-Rinehart
A-algebra. Let us observe
that Lie-Rinehart A-algebras with trivial homomorphism $\alpha: \ele\to $Der(A)
are exactly Lie A-algebras. Therefore the concept of Lie-Rinehart algebras
generalizes the concept of Lie A-algebras. If A=$K$, then Der(A)=0 and there 
is no difference between Lie and Lie-Rinehart algebras. We
denote by ${\ele}{\re}$(A) the category of Lie-Rinehart algebras. One has the
full inclusion
$$ {\ele} {(\mathrm A)} \subset {\ele}{\re}({\mathrm A}),$$
where $\ele({\mathrm A})$ denotes the category of Lie A-algebras. Let us
observe that the kernel of any Lie-Rinehart algebra homomorphism is a Lie
A-algebra.

\begin{Ex} \label{trans} {\rm If ${\frak g}$ is a $K$-Lie algebra acting  
on a commutative $K$-algebra A by derivations (that is there is given a 
homomorphism of Lie $K$-algebras
 $\gamma : {\frak g} \to $ Der(A)), then {\it the
 transformation Lie-Rinehart algebra
of} (${\frak g}$, A) is ${\ele} $ = A $\otimes $ ${\frak g}$ with the Lie
bracket $$[a\otimes g,a'\otimes g'] : = aa' \otimes
[g,g']+a\gamma(g)(a')\otimes g'-a'\gamma(g')(a)\otimes g$$ and with the action
$\alpha :{\ele} \to $Der(A) given by $\alpha(a \otimes g)(a') =
a\gamma(g)(a')$.}
\end{Ex}

\begin{De}{\rm {\it A Lie-Rinehart module} over a Lie-Rinehart A-algebra $\ele$ is a
vector space $M$ together with two operations
$${\ele\t M\to M}, \ \ \ (X,m)\mapsto X(m)$$ and
$$\A\t M\to M\ \ \ (a,m)\mapsto am ,$$ such that the
first one makes $M$ into a module over the Lie $K$-algebra $\ele$ in the sense
of the Lie algebra theory, while the second map makes $M$ into an A-module and
additionally the following compatibility conditions hold
$$(aX)(m)=a(X(m)),$$
$$X(am)=aX(m)+X(a)m.$$
Here $a\in \A$, $m\in M$ and $X\in \ele$.}
\end{De}

It follows that $\A$ is a Lie-Rinehart  module over $\ele$ for any Lie-Rinehart
algebra $\ele$. We let {\bf $({\ele},{\mathrm A})$-mod} be the category of
Lie-Rinehart modules over ${\ele}$.

\subsection{Rinehart cohomology of Lie-Rinehart algebras}
Let $M$ be a Lie-Rinehart module over $\ele.$ Let us recall the definition of
the Rinehart cohomology $H^{\star}_{Rin}({\ele},M)$ of a Lie-Rinehart algebra
${\ele}$ with coefficients in a Lie-Rinehart module $M$ (see \cite{Ri} and
\cite{Hu}). One puts
$$C^n_\A({\ele},M) := Hom_{\mathrm A} (\Lambda^n_{\mathrm A}{\ele},
M),$$ where $\Lambda^*_{\mathrm A}(V)$ denotes the exterior algebra over A
generated by an A-module V. The coboundary map $$\delta : C^{n-1}_\A({\ele},M)
\to C^n_\A({\ele},M)$$ is given by
$$(\delta f)(X_1, \dots,X_n) = (-1)^n \sum_{i=1}^n(-1)^{(i-1)} X_i(f(X_1, 
\dots , \hat{X_i}, \dots,X_n)) + $$
$$ (-1)^n \sum_{j<k} (-1)^{j+k} f([X_i,X_j],X_1, \dots,
\hat{X_i}, \dots,\hat{X_j}, \dots,X_n).$$ Here $X_1, \dots,X_n \in {\ele}, f
\in C^{n-1}_\A({\ele},M)$. By the definition $H^{\star}_{Rin}({\ele},M)$ is the
cohomology of the cochain complex $C^{\star}_\A({\ele},M)$.
 One observes that if $A=K$, then this definition generalizes the classical
 definition of Lie algebra cohomology.
For a general A by forgetting the A-module structure one obtains the canonical
homomorphism
$$H^{\star}_{Rin}({\ele},M)\to H^{\star}_{Lie}({\ele},M),$$
where $H^{\star}_{Lie}({\ele},M)$ denotes the cohomology of $\ele$ considered
as a Lie $K$-algebra. On the other hand if
 A is a smooth commutative algebra, then
$H^{\star}_{Rin}(\mathrm{Der}(\A), \A)$ is isomorphic to the de Rham cohomology
of A (see \cite{Ri} and \cite{Hu}).

 It follows from the definition that one has the following exact sequence
\begin{equation}\label{nolerti}
0\to H^0_{Rin}({\ele},M)\to M\to   \mathrm{Der}_\A({\ele},M)\to
H^1_{Rin}({\ele},M)\to 0,
\end{equation}
where $ {\mathrm{Der}_\A}({\ele},M)$ consists with A-linear maps $d:{\ele}\to
M$ which are derivations from the Lie $K$-algebra ${\ele}$ to $M$. In other
words $d$ must satisfy the following conditions:
$$d(ax)=ad(x), \ a\in \A, x\in \ele,$$
$$d([x,y])=[x,d(y)]-[y,d(x)].$$

For a Lie-Rinehart module $M$ over a Lie-Rinehart algebra $\ele$ one can define
{\it the semi-direct product} ${\ele}\rtimes M$ to be $\ele \oplus M$ as an
A-module with the bracket $[(X,m),(Y,n)]=([X,Y],[X,m]-[Y,n])$.

\begin{Le}\label{kvetader}
Let $\ele$ be a Lie-Rinehart algebra over a commutative algebra {\rm A} and let
$M\in ${\bf $({\ele},{\mathrm A})$-mod}. Then there is a 1-1 correspondence
between the elements of ${\mathrm{Der}_\A}({\ele},M)$ and  the sections (in the
category ${\ele}{\re}$({\rm A})) of the projection $p:{\ele}\rtimes M\to
{\ele}$.
\end{Le}

{\it Proof}.
 Any section $\xi:{\ele}\to {\ele}\rtimes M$ of $p$ has the form 
$\xi(x)=(x,f(x))$ and one easily
shows that $\xi$ is a morphism in ${\ele}{\re}$(A) iff $f\in
{\mathrm{Der}_\A}({\ele},M)$. \qed

\subsection{Abelian and crossed extensions of Lie-Rinehart algebras}

 \begin{De}
Let $\ele$ be a Lie-Rinehart algebra over a commutative algebra ${\mathrm A}$
and let M$\in ${\bf $({\ele},{\mathrm A})$-mod}. An abelian extension of $\ele$
by $M$ is  an exact sequence
$$\xymatrix{0 \ar[r] &{\mathrm M }\ar[r]^i& \ele' \ar[r]^\partial &\ele\ar[r]& 0}$$
where $\ele'$ is a Lie-Rinehart algebra over $\A$ and $\partial$ is a
Lie-Rinehart algebra homomorphism. Moreover, $i$ is an $\A$-linear map and the
following identities hold
$$[i(m),i(n)]=0,$$
$$[i(m),X']=(\partial(X'))m,$$
where $m,n\in $ M and $X'\in \ele'$. An abelian extension is called A-spit if
$\partial$ has an A-linear section.
 \end{De}
 We also need the notion of crossed modules for Lie-Rinehart algebras
introduced in \cite{clpcr}. The following definition is equivalent to one given
in \cite{clpcr}.
\begin{De}
 A crossed module $\partial : {\re} \to {\ele}$ of Lie-Rinehart algebras over
${\mathrm A}$ consists of a Lie-Rinehart algebra ${\ele}$ and a Lie-Rinehart
module $\re$ over $\ele$ together with an $\A$-linear homomorphism $\partial :
{\re} \to {\ele}$ such that for all $r, s \in {\re}, X \in {\ele}, a \in \A$
the following identities hold:
\begin{enumerate}
\item $\partial(X(r)) = [X,\partial(r)]$

\item $(\partial (r))(s)+(\partial(s))(r)=0$

 \item $\partial(r)(a) = 0.$
\end{enumerate}
\end{De}

It follows from this definition that $\re$ is a Lie A-algebra under the bracket
$[r,s]=(\partial(r))(s)$ and $\partial$ is a homomorphism of Lie $K$-algebras.
Moreover, $\Im(\partial)$ is simultaneously a Lie $K$-ideal of ${\ele}$ and an
A-submodule, therefore  $\Coker(\partial)$ is  a  Lie-Rinehart algebra.
Furthermore $\Ker(\partial)$ is an abelian A-ideal of ${\re}$ and the action of
${\ele}$ on $R$ yields a Lie-Rinehart module structure of $\Coker(\partial)$ on
$\Ker(\partial)$.

Let $\pe$ be a Lie-Rinehart algebra  and let M be a Lie-Rinehart module over
${\pe}$. We consider the category $\bf Cross({\pe},{\mathrm M})$, whose objects
are the exact sequences
$$0 \to {\mathrm M} \to {\re} \stackrel{\partial}\to {\ele} 
\stackrel{\upsilon}\to {\pe} \to 0$$
where $\partial : {\re} \to {\ele}$ is a crossed module of Lie-Rinehart
algebras over A and the canonical maps $\Coker( \partial) \to {\pe}$ and
${\mathrm M} \to \Ker( \partial)$ are isomorphisms of Lie-Rinehart algebras and
modules  respectively. The morphisms in the category  $\bf Cross({\pe},{\mathrm
M})$ are commutative diagrams

\begin{center}
\begin{picture}(100,55)
\put(-90,50){$0$}  \put(-40,50){${\mathrm M}$}  \put(10,50){${\re}$}
\put(60,50){${\ele}$} \put(110,50){${\pe}$} \put(160,50){$0$} \put(-90,0){$0$}
\put(-40,0){${\mathrm M}$}  \put(10,0){${\re}'$}  \put(60,0){${\ele}'$}
\put(110,0){${\pe}'$} \put(160,0){$0$}

\put(-80,53){$\vector(1,0){30}$} \put(-25,53){$\vector(1,0){30}$}
\put(25,53){$\vector(1,0){30}$} \put(75,53){$\vector(1,0){30}$}
\put(125,53){$\vector(1,0){30}$} \put(-80,3){$\vector(1,0){30}$}
\put(-25,3){$\vector(1,0){30}$} \put(25,3){$\vector(1,0){30}$}
\put(75,3){$\vector(1,0){30}$} \put(125,3){$\vector(1,0){30}$}

\put(15,43){$\vector(0,-1){30}$} \put(65,43){$\vector(0,-1){30}$}
\put(111,15){$\rule{.5pt}{10mm}$} \put(114,15){$\rule{.5pt}{10mm}$}
\put(-37,15){$\rule{.5pt}{10mm}$} \put(-34,15){$\rule{.5pt}{10mm}$}
\put(5,25){$\alpha$} \put(35,55){${\partial}$} \put(35,5){${\partial}'$}
\put(68,25){$\beta$}
\end{picture}
\end{center}
where $\beta$ is a homomorphism of Lie-Rinehart algebras, $\alpha$ is a
morphism of Lie A-algebras and for any $r \in {\re}, X \in {\ele}$ one has
$$\alpha([X,r]) = [\beta(X), \alpha(r)].$$

Furthermore, we let  ${\bf Cross}_{\mathrm A-spl}({\pe},{\mathrm M})$ be the
subcategory of  $\bf Cross({\pe},{\mathrm M})$ whose objects
and morphisms
split in the category of A-modules,
 in other words one requires  that the epimorphisms
 ${\ele} \to {\pe}$, ${\re} \to \Im( \partial)$,
${\ele}' \to {\pe}'$, ${\re} '\to \Im( \partial)'$,
 $\ele\to \Im(\beta)$,
 $\ele ' \to {\sf Coker}(\beta)$,  $\re\to \Im(\alpha)$,
 $\re ' \to {\sf Coker}(\alpha)$  have
A-linear sections.

\subsection{Main properties of Rinehart cohomologies}

\begin{Th}\label{ziritaditvisebebi} {\rm i)} If $\ele$ is 
projective as an {\rm A}-module, then
$$H^{\star}_{Rin}({\ele},M)\cong \Ext^*_{({\ele},{\mathrm A})-{\bf \mathrm mod}}(A,M).$$

{\rm ii)} If $0\to M_1\to M\to M_2\to 0$ is an exact sequence in the category
{\bf $({\ele},{\mathrm A})$-mod}, then one has a long exact sequence on
cohomology
$$\cdots \to  H^{n}_{Rin}({\ele},M_1) \to H^{n}_{Rin}({\ele},M)\to
 H^{n}_{Rin}({\ele},M_2)\to \cdots$$
provided $0\to M_1\to M\to M_2\to 0$ splits in the category of {\rm A}-modules,
or $\ele$ is projective as {\rm A}-module.

{\rm iii)} The cohomology $H^2_{Rin}({\ele},M)$ classifies such abelian
extensions
$$0\to M\to {\ele} '\to {\ele} \to 0$$
of ${\ele}$ by $M$ in the category of Lie-Rinehart  algebras which splits in
the category of {\rm A}-modules.

 {\rm iv)} For any Lie-Rinehart algebra ${\pe}$ which is projective as an {\rm A}-module
and any Lie-Rinehart module {\rm M} there exists a natural bijection between
the classes of the connected components of the category  ${\bf Cross}_{\mathrm
A-spl}({\pe},{\mathrm M})$
 and $H^3_{Rin}({\pe},{\mathrm M})$.
\end{Th}

{\it Proof.} For the isomorphism of the  part i) see Section 4 of
\cite{Ri}. The part ii) is trivial and for part iii) see
Theorem 2.6 in \cite{Hu}. Finally the part iv), which is in the same spirit as
the classical result for group and Lie algebra cohomology ( see \cite{Lo} and
\cite{KL}), was proved in \cite{clpcr}.

\qed

Let $\frak g$ be a Lie algebra over $K$ and let $M$ be a $\frak g$-module. Then
we have the Chevalley--Eilenberg cochain complex $C^*_{Lie}({\frak g},M)$,
which computes the Lie algebra cohomology (see \cite{CE}):
$$C^n_{Lie}({\frak g},M)=\Hom(\Lambda ^n({\ele}),M).$$
Here $\Lambda ^*$  denotes the exterior algebra defined over $K$.
\begin{Le}\label{transformation} Let ${\frak g}$  be a Lie $K$-algebra acting on a
commutative algebra {\rm A} by derivations and let $\ele$ be the transformation
Lie-Rinehart algebra of $({\frak g},{\mathrm A})$ (see Example \ref{trans}).
Then for any Lie-Rinehart $\ele$-module $M$ one has the canonical isomorphism
of cochain complexes $C^*_\A({\ele},M)\cong C^*_{Lie}({\frak g},M)$ and in
particular the isomorphism
$$H^*_{Rin}({\ele}, M)\cong  H^*_{Lie}({\frak g},M).$$

\end{Le}

{\it Proof}. Since $\ele=$A$\otimes {\frak g} $ one has $Hom_{\mathrm
A}(\Lambda ^n_{\mathrm A}{\ele},M)\cong Hom(\Lambda ^n{\frak g},M) $ and Lemma
follows. \qed

\section{The main construction}
Thanks to Theorem \ref{ziritaditvisebebi} the cohomology theory
$H^*_{Rin}({\ele}, -)$  has good properties only if ${\ele}$ is projective as
an A-module. In this section we introduce the bicomplex $C^{**}({\mathrm A},
\ele, M)$, whose cohomology are good replacement of the Rinehart cohomology
$H^*_{Rin}({\ele}, -)$ for general $\ele$. The idea of the construction is very
simple. First one observes that the transformation Lie-Rinehart algebras (see
Example \ref{trans}) are always free as A-modules, therefore the Rinehart
cohomology of such algebras are the right objects.
Secondly, for any Lie-Rinehart
algebra $\ele$ the two-sided bar construction $B_*(\mathrm{A},
\mathrm{A},{\ele})$ gives rise to a simplicial resolution of $\ele$ in the
category of Lie--Rinehart algebras. Since each term of this resolution is a
transformation Lie-Rinehart algebra one can mix the  Chevalley--Eilenberg
complexes  with the bar resolution to get our bicomplex.

\subsection{A bicomplex for Lie--Rinehart algebras} Let $\ele$ be  a Lie--Rinehart algebra and
let $M$ be a Lie--Rinehart module over $\ele$. We have two cochain complex: the
Rinehart complex $C^*_A({\ele},M)$ and Chevalley-Eilenberg complex
$C^*_{Lie}({\ele},M)$. If one forgets A-module structure on $\ele$, we get a
Lie $K$-algebra acting on A via derivations, thus the construction of Example
\ref{trans} gives a Lie-Rinehart algebra structure on A$\otimes \ele$. One can
iterated this construction to conclude that $\mathrm{A}^{\otimes n}\otimes
\ele$ is also a Lie--Rinehart algebra for any $n\geq 0$. The A-module structure
comes from the first factor, while the bracket is a bit more complicated, for
example for $n=2$, one has
$$[a_1\otimes a_2\otimes X, b_1\otimes b_2\otimes Y]:= a_1b_1\otimes
a_2b_2\otimes [X,Y]+a_1b_1\otimes a_2X(b_2)\otimes Y+$$ $$+a_1a_2X(b_1)\otimes
b_2\otimes Y-a_1b_1\otimes b_2Y(a_2)\otimes X-b_1b_2Y(a_1)\otimes a_2\otimes
X.$$ Let us also recall that the two-sided bar construction
$B_*(\mathrm{A},\mathrm{A},{\ele})$ is a simplicial object, which is
$\mathrm{A}^{\otimes n+1}\otimes {\ele}$ in the dimension $n$, while the face
maps are given by
$$d_i(a_0\otimes \cdots\otimes a_n\otimes X)= a_0\otimes \cdots\otimes
a_ia_{i+1}\otimes\cdots \otimes a_n\otimes X,$$ if $i<n$ and
$$d_n(a_0\otimes \cdots\otimes a_n\otimes X)= a_0\otimes \cdots \otimes a_{n-1}\otimes
\cdots \otimes a_nX,$$ if $i=n$. The degeneracy maps are given by
$$s_i(a_0\otimes \cdots\otimes a_n\otimes X)=a_0\otimes \cdots\otimes a_i\otimes 1\otimes
\cdots a_n\otimes X$$ In fact $B_*(\mathrm{A},\mathrm{A},{\ele})$  is an
augmented simplicial object in the category of Lie-Rinehart algebras, the
augmentation $B_0(\mathrm{A},\mathrm{A},{\ele})=\mathrm{A}\otimes \ele\to
{\ele}$ is given by $(a,X)\mapsto aX$. We can apply the functor
$C^*_\mathrm{A}(-,M)$ on $B_*(\mathrm{A},\mathrm{A},{\ele})$ to get a
cosimplicial object in the category of cochain complexes
$$[n]\mapsto C^*_{\mathrm{A}}(\mathrm{A}^{\otimes n+1}\otimes {\ele},M).$$
Finally we let $C^{**}(\mathrm{A},{\ele},M)$ be the bicomplex associated to
this cosimplicial cochain complex. We let $H^*(\mathrm{ A},{\ele},M)$ be the
cohomology of the corresponding total complex. The augmentation
$B_*(\mathrm{A},\mathrm{A},{\ele})\to {\ele}$ yields the homomorphism
$$\alpha^n:H^*_{Rin}({\ele},M)\to H^*(\mathrm{ A},{\ele},M).$$

The bicomplex $C^{**}(\mathrm{A},{\ele},M)$ has the following alternative
description. According to Lemma \ref{transformation} one has the isomorphism of
complexes:
$$ C^{p*}(\mathrm{A},{\ele},M)\cong C^*_{Lie}(\mathrm{A}^{\t p} \t {\ele},M),$$
where $M$ is considered as a module over $\mathrm{A}^{\t p}\t \ele$ by
$$(a_1\t\cdots \t a_p\t r)m:=(a_1\cdots a_pr)m.$$
To define the horizontal cochain complex structure one observes that elements
of $C^{pq}$ can be identified with functions $f:\mathrm{A}^{\t pq}\t {\ele}^{\t
q}\to M$, which are alternative with appropriate blocks of variables. Then the
corresponding linear map
$$d(f):\mathrm{A}^{\t (p+1)q}\t {\ele}^{\t q}\to M$$
is  given by
$$df(a_{01},\cdots ,a_{0q},a_{11},\cdots, a_{1q},\cdots a_{p1},\cdots a_{pq},
X_1,\cdots,X_q)=$$
$$a_{01}\cdots a_{0q}f(a_{11},\cdots, a_{1q},\cdots a_{p1},\cdots a_{pq},
X_1,\cdots,X_q)+$$
$$+\sum_{0\leq i<p}(-1)^{i+1}f(a_{01},\cdots ,a_{0q},\cdots,a_{i1}a_{i+1,1}\cdots,a_{iq}a_{i+1,q},\cdots a_{p1},\cdots a_{pq},X_1,\cdots,X_q)+$$
$$(-1)^{p+1}f(a_{01},\cdots ,a_{0q},\cdots,a_{p-1,1},\cdots ,a_{p-1,q},a_{p1}X_1,\cdots,a_{pq}X_q).$$

\begin{Th}\label{zlivs} {\rm i)} The homomorphism
$$\alpha^n:H^n_{Rin}({\ele},M)\to H^n(\mathrm{A},{\ele},M)$$
is isomorphism for $n=0,1$.  The homomorphism $\alpha^2$ is a monomorphism.
Moreover $\alpha ^n$ is an isomorphism for all $n\geq 0$ provided $\ele$ is
projective over $\mathrm{A}$.

{\rm ii)} If $0\to M_1\to M\to M_2\to 0$ is an exact sequence in the category
{\bf $({\ele},{\mathrm A})$-mod}, then one has a long exact sequence on
cohomology
$$\cdots \to  H^{n}(\mathrm{A},{\ele},M_1) \to H^{n}(\mathrm{A},{\ele},M)\to
 H^{n}(\mathrm{A}, {\ele},M_2)\to \cdots.$$

{\rm iii)} The cohomology $H^2(\mathrm{A},{\ele},M)$ classifies all abelian
extensions
$$0\to M\to {\ele} '\to {\ele} \to 0$$
of ${\ele}$ by $M$ in the category of Lie-Rinehart  algebras.

 {\rm iv)} For any Lie-Rinehart algebra ${\ele}$
and any Lie-Rinehart module {\rm M} there exists a natural bijection between
the classes of the connected components of the category  ${\bf
Cross}({\ele},{\mathrm M})$
 and $H^3(\mathrm{A},{\ele},{\mathrm M})$.
\end{Th}

{\it Proof}. i) The statement is obvious for $n=0,1$. For $n=2$ it follows from the 
part iii) below and part iii) of Theorem \ref{ziritaditvisebebi}. It remains to 
prove the last assertion. It is well-know that the augmentation 
$B_*(\mathrm{A},\mathrm{A},{\ele})\to \ele$ is a homotopy equivalence in the 
category of simplicial vector spaces, thanks to the existence of the extra 
degeneracy map given by
 $s(a_0\otimes \cdots\otimes a_n\otimes X)=1\otimes a_0\otimes 
\cdots\otimes a_n\otimes
X$. However $s$ is not A-linear and therefore in general 
$B_*(\mathrm{A},\mathrm{A},{\ele})\to \ele$ is only a weak 
equivalence in the category of  simplicial A-modules. 
Assume now $\ele$ is
projective as an A-module, then  $B_*(\mathrm{A},\mathrm{A},{\ele})\to \ele$ is 
a homotopy equivalence  in the category of  simplicial A-modules and therefore, 
for each $k\geq 0$ the induced map $\Lambda^k_\mathrm{A}(B_*(\mathrm{A},
\mathrm{A},{\ele}))\to
\Lambda^k_\mathrm{A}({\ele})$ is  a homotopy equivalence  
in the category of  simplicial A-modules, which implies that 
the same is true after applying
the functor ${\sf Hom}_
\mathrm{A}(-,M)$. Thus for each $k\geq 0$ the induced map
$C^k_\mathrm{A}({\ele},M)\to C^k_\mathrm{A}
(B_*(\mathrm{A},\mathrm{A},{\ele}))$ is a weak equivalence of cosimplicial objects and the
comparison theorem for bicomplexes yields the result.

ii) Since Hom and exterior powers involved  in $C^m_{Lie}({\frak g},M)$ are taken over $K$ it 
follows that for each $p$ and $q$  the functor $C^q_{Lie}(\mathrm{A}^p\otimes {\ele},-)$ is
exact and the result follows.

iii)  Thanks to a well-known  fact from topology we can use the normalized (in 
the simplicial direction) cochains to compute $H^*(\mathrm{A},{\ele},M)$. Having 
this in mind we
have $H^2(\mathrm{A},{\ele},M)=Z^2/B^2$, 
where $Z^2$  consists with pairs $(f,g)$ such that $f:\Lambda ^2({\ele})\to M$ is a 
Lie 2-cocycle and $g:\mathrm{A}\t {\ele}
\to M$ is a linear map such that $g(1,X)=0$, $$ag(b,X)-g(ab,X)+g(a,bX)=0$$ and 
$$abf(X,Y)-f(aX,bY)=$$
$$aXg(b,Y)-bYg(a,X)-g(ab,[X,Y])-g(aX(b),Y)+g(bY(a),X).$$
Here $a,b\in{\mathrm A}$ and $X,Y\in \ele$. Moreover $(f,g)$ belongs to
$B^2$ iff there exist a linear map $h:{\ele}\to M$ such that
$f(X,Y)=Xh(Y)-h([X,Y])-Yh(X)$ and $g(a,X)=ah(X)-h(aX).$
 Starting with $(f,g)\in Z^2$ we construct an abelian extension of $\ele$ by $M$ by
putting ${\pe}=M\oplus \ele$ as a vector space. An A-module structure on ${\pe}$ is
given by $a(m,X)=(am+g(a,X),aX)$, while a Lie bracket on $\pe$ is given by
$[(m,X),(n,Y)]=(Xn-Ym+f(X,Y),[X,Y])$. Conversely, given an abelian extension $(\pe)$ and
a $K$-linear section $h:{\ele}\to \pe$ we put
$f(X,Y):=[h(X),h(Y)]-h([X,Y])$ and $g(a,X):=h(aX)-ah(X)$. One easily checks that
$(f,g)\in  Z^2$ and one gets iii). 

iv) Similarly, we have
$H^3(\mathrm{A},{\ele},M)=Z^3/B^3$. Here $Z^3$  consists with triples
$(f,g,h)$ such that $f:\Lambda ^3 ({\ele})\to M$ is a Lie 3-cocycle,
 $g:\Lambda ^2(\mathrm{A}\t {\ele})\to M$ and $h:\mathrm{A}
\t \mathrm{A}\t {\ele}\to M$ are linear maps and the following relations hold:
$$f(aX,bY,cZ)- abcf(X,Y,Z)= $$
$$=aXg(b,c,Y,Z)-bYg(a,c,X,Z)+cZg(a,b,X,Y)-g(ab,c,[X,Y],Z)+$$
$$+g(aX(b),c,Y,Z)- g(bY(a),c,X,Z)+g(ac,b,[X,Y],Y)-g(aX(c),b,Z,Y)+$$
$$+g(cZ(a),b,X,Y)-g(bc,a,[Y,Z],X)+g(bY(c),a,Z,X)-g(cZ(b),a,Y,X)$$
and
$$abXh(c,d,Y)-cdYh(a,b,X)-h(ac,hd,[X,Y])-h(ac,bX(d),Y)-$$
$$-h(abX(c),d,Y)+h(ac,dY(b),X)-h(cdY(a),b,X)=$$
$$=abg(c,d,X,Y)-g(ac,bd,X,Y)+g(a,b,cX,dY).$$
Moreover, $(f,g,h)$ belongs to $B^3$ iff there exist a linear maps
$m: \Lambda ^2({\ele})\to M$ and $n:\mathrm{A}\t {\ele}\to M$ such that
$$f(X,Y,Z)=Xm(Y,Z)-Ym(X,Z)+Zm(X,Y)-$$
$$-m([X,Y],Z)+m([X,Z],Y)-m([Y,Z],X),$$
$$g(a,b,r,s)=abm(X,Y)-m(aX,bY)-aXn(b,Y)+$$
$$bYn(a,X)+n(ab,[X,Y])+n(aX(b),Y)+n(bY(a),X)$$
and
$$h(a,b,X)=an(b,X)-n(ab,X)+n(a,bX).$$

Let $$0\to M\to {\re} \buildrel \partial \over \to {\pe} \buildrel \pi \over \to
{\ele}\to 0$$ be a crossed extension. We put $V:=\Im(\partial)$ and consider
$K$-linear sections $p:{\ele}\to \pe$ and $q:V\to \re$ of $\pi: \pe\to \ele$ and
$\partial:{\re }\to V$ respectively. Now we define $t:{\ele}\t {\ele}\to \re$
and $s:\mathrm{A}\t {\ele}\to \re$.
by $t(X,Y):=q([p(X),p(Y)])-p([X,Y])$ and $s(a,X):=q(ap(X)-p(aX))$. 
Finally we define three functions as follows.
The function $f: \Lambda ^3({\ele})\to M$ is given by
$$f(X,Y,Z):=p(X)g(Y,Z)-p(Y)g(X,Z)+p(Z)g(X,Y)-$$
$$-g([X,Y],Z)+g([X,Z],Y)-g([Y,Z],X).$$
The function $g:\Lambda ^2(\mathrm{A}\t {\ele})\to M$
is given by $$g(a,b,X,Y):= p(aX)s(b,Y)-p(bY)s(a,X)-p(ab,[X,Y])-
$$
$$-p(aX(b),Y)+p(bY(a),X)-
t(aX,bY)+abt(X,Y),$$
while the function 
and $h:\mathrm{A}\t \mathrm{A}\t {\ele}\to M$ is given by
$$h(a,b,X):= as(b,X)-s(ab,X)+s(a,bX).$$
Then $(f,g,h)\in Z^3$ and the corresponding class in $\H^3(A,R,M)$
depends only on the connected component of a given crossed extension. Thus we
obtain a well-defined map $ {\bf Cros}(A,R,M)\to \H^3(A,R,M)$ and a standard
argument (see \cite{KL})  shows that it is an isomorphism.

\section{Triple cohomology of Lie--Rinehart algebras}
In this section we prove that the chomology theory developped in the previous 
section in canonically isomorphic to the triple cohomology of Barr-Beck \cite{bb} applyed
to Lie-Rinehart algebras.

\subsection{Cotriples and cotriple resolutions} The general notion of
(co)triples (or (co)monads, or (co)standard construction) and (co)triple resolutions due to  Godement \cite{godem} and then was developed in \cite{bb}. 
Let ${\cal C}$ be a category.
A {\it cotriple} on ${\cal C}$ is an endofunctor $T:{\cal C}\to {\cal C}$
together with natural transformations $\epsilon: T\to 1_{\cal C}$ and $\delta: T\to T^2$
satisfying the counit and the coassociativity properties. Here $T^2=T\circ T$ and a
similar meaning has $T^n$ for all $n\geq 0$. For example, assume $U:{\cal C}\to {\cal B}$
is a functor which has a left adjoint functor $F:{\cal B}\to {\cal C}$.
Then there is a cotriple structure on $T=FU:{\cal C}\to {\cal C}$ such that
$\epsilon$ is the counit of the adjunction.
Given a cotriple $T$ and an object $C$ one can associate a simplicial object $T_*C$ in the category ${\cal C}$, known as {\it Godement or cotriple resolution of $C$}.
Let us recall that $T_nC=T^{n+1}C$ and the face and  degeneracy operators are given respectively by $\partial_i=T^i\epsilon T^{n-i}$ and
$ s_i= T^i\delta T^{n-i}$. To explain why it is called resolution, consider the case when $T=FU$ is associated to the pair of adjoint functors. Then
 firstly $\epsilon $ yields a morphism  $T_*C\to C$ from the simplicial object $T_*C$ to the constant simplicial object
$C$ and secondly the induced morphism $U(T_*C)\to U(C)$ is a homotopy equivalence in the category of simplicial objects in $\cal B$.
The cotriple cohomology is now defined as follows. Let $M$ be an abelian group object in the category ${\cal C}/C$ of arrows $X\to C$ then
$Hom_{{\cal C}/C}(T_*C,M)$ is a cosimplicial abelian group, which can be seen also as a cochain complex. Thus $H^*(Hom_{{\cal C}/C}(T_*C,M))$ is a meaningful and they are denoted  by $H^*_T(C,M)$.
Of the special interest is the case, when $T=FU$ is associated to the pair of adjoint functors and the functor $U:{\cal C}\to {\cal B}$ is {\it tripliable} \cite{bb}. In this case the
category $\cal C$ is completely determined by the triple $E=UF:{\cal B}\to  {\cal B}$.
Because of this fact $H^*_T(C,M)$ in this case  are known as {\it triple cohomology
of $C$ with coefficients in $M$}.

\subsection{Free Lie--Rinehart Algebras}
We wish to apply these general constructions to Lie--Rinehart algebras. One has
the functor
$$U: {\ele}{\re}({\rm A})\to {\rm Vect/ Der(A)}$$
which assigns $\alpha:{\ele}\to {\rm Der(A)}$ to a Lie--Rinehart algebra
${\ele}$. Here Vect/ Der(A)  is the category of $K$-linear maps $\psi:V\to $
Der(A), where $V$ is a vector space over $K$. A morphism $\psi\to \psi_1 $ in
Vect/Der(A)  is a $K$-linear map $f:V\to V_1$ such that $\psi=\psi_1\circ f$.
Now we construct the functor
$$F:{\rm Vect/Der(A)}\to {\ele}{\re}({\rm A})$$
as follows. Let $\psi:V\to $Der(A)  be a $K$-linear map. We let ${\bf L}(V)$ be
the free Lie $K$-algebra generated by $V$. Then one has the unique Lie
$K$-algebra homomorphism ${\bf L}(V)\to \ $Der(A) which extends the map $\psi$,
which is still denoted by $\psi$. Now we can apply the construction from
Example \ref{trans} to get a Lie--Rinehart algebra structure on A$\otimes {\bf
L}(V)$. We let $F(\psi)$ be this particular Lie--Rinehart algebra and we call
it {\it the free Lie--Rinehart algebra generated by $\psi$}. In this way we
obtain the functor $F$, which is the left adjoint to $U$.

\begin{Le}\label{free} Let ${\ele}$ be a free Lie--Rinehart 
algebra generated by $\psi:V\to $Der(A) and let $M$ be any Lie--Rinehart 
module over ${\ele}$. Then
$$H^i_{Rin}({\ele},M)=0, \ \ i>1.$$
\end{Le}

{ \it Proof}.  By our construction  $\ele$ is a transformation Lie--Rinehart
algebra of $({\bf L}(V),{\rm A})$. Thus one can apply Lemma
\ref{transformation} to get an isomorphism $H^*_{Rin}({\ele},M)\cong
H^*_{Lie}({\bf L}(V),M)$ and then one can use the well-known  vanishing result
for free Lie algebras. $\Box$

\subsection{The cohomology  $H^*_{LR}({\ele},M)$}
Since we have a pair of adjoint functors we can take the composite
$$T=FU:{\ele}{\re}({\rm A})\to {\ele}{\re}({\rm A})$$ which is a cotriple. Thus
for any Lie--Rinehart algebra $\ele$ one can take the cotriple resolution
$T_*({\ele})\to {\ele}$. It follows from the construction of the cotriple
resolution that each component of $T_*({\ele})$ is a free Lie--Rinehart
algebra. Moreover, according to the general properties of the cotriple
resolutions  the natural augmentation $T_*({\ele})\to {\ele}$ is a homotopy
equivalence in the category of simplicial vector spaces. It follows that
$T_*({\ele})\to {\ele}$ is a weak homotopy equivalence in the category of
A-modules.

Let $M$ be a $\ele$-module. Then $M$ is also a module over $T_n({\ele})$ for any $n\geq 0$
thanks to the augmentation morphism $T_*({\ele})\to {\ele}$. Thus one can form the following 
bicomplex $$C^*_{\mathrm{A}}(T_*({\ele}),M)$$
which is formed by the degreewise applying the Rinehart cochain complex. 
The cohomology of the total complex of the bicomplex $C^*_{\mathrm{A}}(T_*({\ele}),M)$
is denoted by $H^*_{LR}({\ele},M)$. 


\begin{Le}\label{firstiso} For any  Lie--Rinehart algebra ${\ele}$ and any Lie--Rinehart module
$M$ one has a natural  isomorphism
$$H^*(\mathrm{A},{\ele},M)\cong H^*_{LR}({\ele},M).$$
\end{Le}

{\it Proof}. We denote by $C^*(\mathrm{A},{\ele},M)$ the total complex associated to the bicomplex
$C^{**}(\mathrm{A},{\ele},M)$. Recall that it comes with a natural cochain 
map $$C^*_{\mathrm{A}}({\ele},M)\to C^*(\mathrm{A},{\ele},M)$$
which is quasi-isomorphism provided $\ele$ is projective as an A-module. Let us apply 
$C^*(\mathrm{A},-,M) $ on $T_*({\ele})$ degreewise. Then we obtain the morphism of bicomplex
$$C^*_{\mathrm{A}}(T_*({\ele}),M)\to C^*(\mathrm{A}, T_*({\ele}),M)$$
which is quasi-isomorphism because each $T_n({\ele})$ is free as an A-module.
It remains to show that the augmentation $T_*({\ele})\to \ele$ yields the quasi-isomorphism
$$C^*(\mathrm{A},{\ele},M)\to  C^*(\mathrm{A}, T_*({\ele}),M).$$
To this end, one observes that $T_*({\ele})\to \ele$ is a quasi-isomorphism thanks to the
general properties of cotriple resolutions and therefore is a homotopy equivalence in the category
of simplicial vector spaces. Thus the same is true for $\Lambda^n(T_*({\ele}))\to \Lambda ^n({\ele})$
and therefore $C^n(\mathrm{A},{\ele},M)\to C^n(\mathrm{A}, T_*({\ele}),M)$ is also 
a homotopy-equvalence for each $n$ and the result follows from the comparison theorem of bicomplexes.

\subsection{Triple cohomology and $H^*_{LR}({\ele},M)$}
According to the Backs triplibility criterion 
the functor $U: {\ele}{\re}({\rm A})\to {\rm Vect/ Der(A)}$ is
tripliable, so we have also the triple cohomology theory for Lie--Rinehart
algebras. Let $\ele$ be a Lie--Rinehart algebra. There is an equivalence  from
the category of Lie--Rinehart modules over $\ele$ to the category of abelian
group objects in ${\ele}{\re}({\rm A})/{\ele}$, which assigns the projection
${\ele}\rtimes M\to \ele$ to $M\in ({\ele},{\rm A})-{\bf \rm mod}$. 
Having in
mind this equivalence, Lemma \ref{kvetader} says that for any object ${\pe} \to
{\ele}$ of ${\ele}{\re}({\rm A})/{\ele}$ the homomorphisms from ${\pe} \to
{\ele}$ to ${\ele}\rtimes M\to \ele$ in the category of abelian group objects
in ${\ele}{\re}({\rm A})/{\ele}$ is nothing else, but ${\rm Der_A}({\pe},M)$.
Therefore the triple cohomology $H_T^*({\ele},M)$ is the same as $H^q({\rm
Der_A}(T_*({\ele}),M))$.

\begin{Th}\label{seciso}
For any Lie--Rinehart algebra $\ele$ and any $\ele$-module $M$ there is a
natural isomorphism:
$$H^{q+1}_{LR}({\ele},M)\cong H_T^q({\ele},M), \ \ q>0.$$
In other words the cotriple cohomology of $\ele$ with coefficients in 
$M$ is isomorphic to the cohomology
 $H_{LR}^*({\ele},M)$ up to shift in the dimension.
\end{Th}

{\it Proof}.
As usual with bicomplex we have a spectral sequence
$$E^2_{pq}\Longrightarrow H^*_{LR}({\ele},M)$$
where $E^2_{pq}$ is obtained in two steps: First one takes $p$-th homology 
in each $C^*(T_q({\ele}),M)$, $q\geq 0$ and then one takes the
 $q$-th homology.
But $C^*(T_q({\ele}),M)$ is just the Rinehart complex of $T_q({\ele})$. Since $T_q({\ele})$
is free we can use  Lemma \ref{free} to conclude that
$E^1_{pq}=0$ for all $p\geq 2$. According to the exact
sequence (\ref{nolerti}) one has also an exact sequence
$$0\to E^1_{0q}\to M\to {\rm Der_A}(T_q({\ele}),M)\to E^1_{1q}\to 0$$
One observes that $ E^1_{0*}$ and $M$ are  constant cosimplicial vector spaces 
and therefore $E^2_{0q}=0$ for all $q>0$. Thus we get
$$H^{q+1}_{LR}({\ele},M)\cong E^2_{1q}\cong H^q({\rm Der_A}(T_*({\ele}),M)), \ q>0.$$

\section{The canonical class of associative algebras}
Let $S$ an associative algebra over $K$. We let $A$ be the center of $S$.
 As an application of our results we construct a canonical class $o(S)\in
H^3(\mathrm{A},H^1(S,S),\mathrm{A})$, where $S$ is an associative algebra and $H^*(S,S)$
denotes the Hochschild cohomology of $S$.

Let us first recall the definitions of the zero and the first dimensional Hochschild cohomology
involved in this construction.
Let $S$ be an associative $K$-algebra. {\it A $K$-derivation} $D: S\to S$ is a $K$-linear
map, such that $D(ab)=D(a)b+aD(b)$. We let
Der($S$) be the set of all $K$-derivations. It has a natural Lie $K$-algebra 
structure, where the bracket is defined via the commutator $[D,D_1]=DD_1-D_1D$.
There is a canonical $K$-linear map
 $${\rm ad}:S \to {\rm Der}(S)$$
given by ${\rm ad}(s)(x)=sx-xs$, $s,x\in S$. Then the zero and the first dimensional Hochschild
cohomology groups are defined via the exact sequence:
\begin{equation}\label{canon}
0\to H^0(S,S) \to S \buildrel {\rm ad} \over
\longrightarrow {\rm Der}(S)\to H^1(S,S)\to 0
\end{equation}
It follows that  A$=H^0(S,S)$ is the center of $S$. We claim that Der($S$) is a
Lie--Rinehart algebra over A. Indeed, the action of A is defined by
$(aD)(s)=aD(s)$, $D\in $ Der($S$), $s\in S$, $a\in $ A, while the homomorphism
$\alpha:$ Der($S$)$\to $Der(A) is just  the restriction. To see that $\alpha$
is well-defined, it suffices to show that $D(A)\subset A$ for any $D\in $
Der($S$). To this end, let us observe that for any $s\in S$ and $a\in A$ one
has
$$D(a)s-sD(a)=(D(as)-aD(s))-(D(sa)-D(s)a)=0$$
and therefore $D(a)\in A$. On the other hand  the commutator $[s,t]=st-ts$
defines a Lie $A$-algebra structure on $S$ and ${\rm ad}:S \to {\rm Der}(S)$ is
a Lie $K$-algebra homomorphism. Actually more is true: ${\rm ad}$ is a crossed
module of Lie--Rinehart algebras over A, where the action of the Lie--Rinehart
algebra Der($S$) on $S$ is given by $(D,s)\mapsto D(s)$. It follows that
$H^1(S,S)=\Coker(ad:S\to {\rm Der}(S))$ is also a Lie--Rinehart algebra over A
and A=$\Ker(ad:S\to {\rm Der}(S))$ is a Lie--Rinehart module over $H^1(S,S)$.
In particular the groups $ H^*(\mathrm{A},H^1(S,S),\mathrm{A})$ 
are well-defined. According to
Theorem \ref{zlivs} the vector space $ H^3(\mathrm{A},H^1(S,S),\mathrm{A})$ classifies the
crossed extension of $H^1(S,S)$ by $A$. By our construction the exact sequence
(\ref{canon}) is one of such extension and therefore it defines a canonical
class $o(S)\in H^3(\mathrm{A},H^1(S,S),\mathrm{A})$.
\begin{Le} $o(S)$ is a Morita invariant.
\end{Le}
{\it Proof}. Let $R$ be the $K$-algebra of $n\times n$ matrices. We have
to prove that $o(S)=o(R)$. Let $D$ be a derivation of $S$. We let $g(D)$ be
the derivation of $R$ which is componentwise extension of $D$. Furthermore,
for an element $s\in S$ we let $f(s)$ be the diagonal matrix with $s$ on diagonals. 
Then one has the following commutative diagram
$$\xymatrix{S\ar[d]_f\ar[r]^{\rm ad}&{\rm Der}(S)\ar[d]^g\\
R\ar[r]^{\rm ad}&{\rm Der}(R)}$$
in the category ${\ele}{\re}$(A) and the result follows from the fact that Hochschild cohomology is a Morita invariant.$\Box$

Let us observe that if $S$ is a smooth commutative algebra, then $A=S$ and
$H^3(\mathrm{A},H^1(S,S),\mathrm{A})$ is isomorphic to the de Rham cohomology of $S$ (of
course $o(S)=0$  in this case). So, in general one can  consider the groups
$H^3(\mathrm{A},H^1(S,S),\mathrm{A})$ as a sort of noncommutative de Rham cohomology.

By forgetting $A$-module structure, one obtains an element
 $$o'(S)\in H^3_{Lie}(H^1(S,S),A).$$
These groups and probably the corresponding elements
can be compute in many cases using the results of Strametz  \cite{St}.

\bigskip
\centerline{\bf Acknowledgments}

\bigskip
The authors were supported by MCYT, Grant BSFM2003-04686-C02.
The third author is very grateful to Universities of Santiago and
Vigo for hospitality.
The two first authors were supported also by DGES BFM2000-0523.

\begin{center}

\end{center}

\end{document}